\documentclass{article}

\usepackage[doc]{optional}
\usepackage{color}
\usepackage{float}
\usepackage{soul}
\usepackage{graphicx}
\definecolor{labelkey}{rgb}{0,0.08,0.45}
\definecolor{refkey}{rgb}{0,0.6,0.0}
\definecolor{Brown}{rgb}{0.45,0.0,0.05}
\definecolor{lime}{rgb}{0.00,0.8,0.0}
\definecolor{lblue}{rgb}{0.5,0.5,0.99}
\usepackage{amsmath,amssymb,amsthm, fancyhdr,enumerate,algorithm,algpseudocode,multicol,tikz}
\usepackage[margin=1.1in]{geometry}
\usepackage{hyperref}
\usepackage{cite}
\usepackage[shortlabels]{enumitem}

\definecolor{labelkey}{rgb}{0.6,0.6,0.6}
\definecolor{refkey}{rgb}{0,0.6,0.0}

\newtheorem{theorem}{Theorem}[section]

\theoremstyle{definition}

\newtheorem{example}[theorem]{Example}

\theoremstyle{remark}
\newtheorem{remark}[theorem]{Remark}

\def\disp{\displaystyle}

\def\tto{\;{\lower 1pt \hbox{$\rightarrow$}}\kern -10pt
	\hbox{\raise 2pt \hbox{$\rightarrow$}}\;}

\def\epsilon{\varepsilon}

\def\R{\Bbb R}

\begin{document}
	
	\title{Minimax Theorems for Possibly Nonconvex Functions}
	
	\author{
		Nguyen Nang Thieu\thanks{Institute of Mathematics, Vietnam Academy of Science and Technology, 18 Hoang Quoc Viet, Hanoi 10307, Vietnam. Email: \texttt{nnthieu@math.ac.vn}} \and
		Nguyen Dong Yen\thanks{Institute of Mathematics, Vietnam Academy of Science and Technology, 18 Hoang Quoc Viet, Hanoi 10307, Vietnam. Email: \texttt{ndyen@math.ac.vn}}
	}
	
	\date{} 
	\maketitle


\begin{abstract} This paper establishes three minimax theorems for possibly nonconvex functions on Euclidean spaces or on infinite-dimensional Hilbert spaces. The theorems also guarantee the existence of saddle points. As a by-product, a complete solution to an interesting open problem related to continuously differentiable functions is obtained. The obtained results are analyzed via a concrete example. 
\end{abstract}



\section{Introduction}\label{Sect-1}

Originating from John von Neumann's celebrated work~\cite{von Neumann_1928} on game theory, \textit{minimax theorems} are widely recognized as important tools in mathematics and mathematical applications. Note that von Neumann~\cite[Sects.~2 and~3]{von Neumann_1928} established a min-max theorem for zero-sum two-person games. Two proofs of the theorem can be found in the book by Barron~\cite[Subsect.~1.2.1]{Barron_book_2013}. From the result, it follows that any matrix two-person game with mixed strategies possesses solutions; see, e.g.,~\cite[Theorem~1.3.4]{Barron_book_2013}. 
	
\medskip
Roughly speaking, a minimax theorem gives some conditions assuring the validity of the equality \begin{equation}\label{minimax_eq} \disp\sup_{y\in Y}\,\disp\inf_{x\in X}f(x,y)=\disp\inf_{x\in X}\,\disp\sup_{y\in Y}f(x,y),\end{equation} where $f:X\times Y\to\mathbb R$ is a function of two variables and $X, Y$ are some sets. It is well known that, without additional assumptions, the following inequality is true \begin{equation}\label{minimax_ine} \disp\sup_{y\in Y}\,\disp\inf_{x\in X}f(x,y)\leq\disp\inf_{x\in X}\,\disp\sup_{y\in Y}f(x,y);\end{equation} so the quantity $$\alpha:=\disp\inf_{x\in X}\,\disp\sup_{y\in Y}f(x,y)-\disp\sup_{y\in Y}\,\disp\inf_{x\in X}f(x,y)$$ is always nonnegative. One may have $\alpha>0$ or even $\alpha=+\infty$. Note that the inequality~\eqref{minimax_ine} can easily be proved. 

\medskip
Certain minimax theorems (for instance, von Neumann's minimax theorem) may include a conclusion about the existence of a \textit{saddle point}, that is a point $(\bar x,\bar y)\in X\times Y$ such that 
	\begin{equation}\label{saddle_point} f(\bar x,y)\leq f(\bar x,\bar y)\leq f(x,\bar y)\quad\; \forall x\in X,\ \forall y\in Y.\end{equation} It is not difficult to show that the existence of a saddle point assures that	
	\begin{equation}\label{minimax_eq1c} \disp\sup_{y\in Y}\,\disp\inf_{x\in X}f(x,y)=f(\bar x,\bar y)=\disp\inf_{x\in X}\,\disp\sup_{y\in Y}f(x,y);\end{equation} hence the fulfillment of~\eqref{saddle_point} for a point $(\bar x,\bar y)\in X\times Y$ implies~\eqref{minimax_eq}. The converse is not true in general. If~\eqref{minimax_eq1c} (resp.,~\eqref {minimax_eq}) holds, then the number~$f(\bar x,\bar y)$ (resp., the number~$\disp\sup_{y\in Y}\,\disp\inf_{x\in X}f(x,y)$, which is equal to $\disp\inf_{x\in X}\,\disp\sup_{y\in Y}f(x,y)$) is called the \textit{saddle value} of $f$; see, e.g.,~\cite[p.~196]{Aubin_1979}.

  \medskip
  Among numerous applications of minimax theorems of different types, we may refer to the existence of solutions in two-person games~(see, e.g., Aubin~\cite[Chapter~7]{Aubin_1979}, Aubin and Ekeland~\cite[Sect.~2 in Chapter~6]{Aubin_Ekeland_1984}, Barron~\cite[Chapters~1 and~2]{Barron_book_2013}), the Lipschitz continuity and a fundamental property of the effective domains of generalized polyhedral convex multifunctions (see~\cite[Theorem~2.207]{BS_2000} and~\cite[Proof of Theorem~3.1]{LTY_2024}), properties of maximal operators in reflexive Banach spaces~\cite{Simons_1998},  differential stability of parametric optimization problems~(see~\cite{Aubin_1980} and~\cite[Theorem~3.2]{DY_1991}), as well as the solution existence of generalized quasi-variational inequalities in infinite-dimensional normal spaces~(see~\cite[Theorem~2.1]{CY_1997}). 
  	
  \medskip
  Observe that the concept of minimax variational inequality, which can serve as a good tool for studying minimax problems given by convex sets and differentiable functions, was suggested in~\cite{HY_2011}.

\medskip
Typical assumptions of a minimax theorem with an assertion of the type~\eqref{minimax_eq} involving a function $f:C\times D\to\mathbb R$, where $C$ and $D$ are convex subsets of topological vector spaces $X$ and $Y$, respectively, are as follows:
	\begin{itemize}
		\item[(a1)] For each $y\in D$, the function $f(.,y):C\to\mathbb R$ is convex;
		\item[(a2)] For each $x\in C$, the function $f(x,.):D\to\mathbb R$ is concave
	\end{itemize}
	(see, for instance,~\cite[Theorem~3.1]{Simons_1998}). Sometimes, the convexity in~(a1) is relaxed to the quasiconvexity (or pseudoconvexity) and the concavity in~(a2) is weakened to the quasiconcavity (or pseudoconcavity). As an example, we may refer to Sion's minimax theorem~(see~\cite[Theorem~3.4]{Sion_1958} and also~\cite[Theorem~7 in Chapter~7]{Aubin_1979}). 

\medskip
	When no assumptions on partial convexity (or generalized convexity) and partial concavity (or generalized concavity) are made on the function under consideration, we will have \textit{minimax theorems for  possibly nonconvex functions}. Most of the topological minimax theorems (see, e.g., Tuy~\cite{Tuy_1974,Tuy_2004,Tuy_2008,Tuy_2011} and Ricceri~\cite{Ricceri_1993,Ricceri_1998}) are of this type.
	
\medskip
	As indicated in the title of this paper, herein we are interested in obtaining minimax theorems for  possibly nonconvex functions. The results arose from our great efforts to solve an open problem. Namely, in the year 2024, Professor Biagio Ricceri~\cite{Ricceri_2024b} sent us via emails three very interesting open problems related to continuously differentiable functions ($C^1$~functions for brevity) and continuous functions. The first problem is formulated as follows.
	
\medskip
\noindent \textbf{Problem 1:} \textit{Let $n\geq 2$. Find a $C^1$ function $J:\R^n \to \R$ and a convex set $Y \subset \R^n$ satisfying the following conditions:
	\begin{itemize}
		\item [\rm (a)] $\lim\limits_{\|x\|\to+\infty} J(x) =+\infty$;
		\item [\rm (b)] $\nabla J$ is Lipschitzian, with Lipschitz constant $L$;
		\item[\rm (c)] one has
		$$\inf_{x\in \R^n} J(x) -\dfrac{L}{2}\inf_{y\in Y}\|y\|^2< \inf_{x\in \R^n}\, \sup_{y\in Y} \left(J(x+y)-\dfrac{L}{2}\|y\|^2\right).$$
	\end{itemize}
} 

The aim of the present paper is to provide a complete solution to the problem. It turns out that if the conditions in~(a) and~(b) are fulfilled, then the inequality in~(c) must be violated. More precisely, instead of the strict inequality, we will have an equality like~\eqref{minimax_eq} with a suitable defined function $f(x,y)$. Our study of Problem~1 yields three minimax theorems for possibly nonconvex functions on $\mathbb R^n$ or on infinite-dimensional Hilbert spaces, which also include assertions about the existence of saddle points. We will use some tools from optimization theory to prove these minimax theorems. To analyze the obtained results, we will present a concrete example. 

\medskip
Our main results are established in Section~\ref{Sect-2}. Section~\ref{Sect-3} is devoted to an illustrative example. Some concluding remarks are given in Section~\ref{Sect-4}.

\section{Main Results}\label{Sect-2}

The following remarks help us to have a closer look at Problem~1. 

\begin{remark}
\rm{Suppose that we can find a function $J$ and a set $Y$ that solve Problem~1. If $Y$ is not closed, we can replace it with its closure, $\overline{Y}$, while still ensuring that $J$ and $\overline{Y}$ solve Problem~1. This is justified by the following assertions:
	\begin{itemize}
		\item[\rm (i)] $\inf\limits_{y\in Y} \|y\|^2 =\inf\limits_{y\in \overline{Y}} \|y\|^2$.
		\item[\rm (ii)] $\sup\limits _{y\in Y} \left(J(x+y)-\dfrac{L}{2}\|y\|^2\right)=\sup\limits _{y\in \overline{Y}} \left(J(x+y)-\dfrac{L}{2}\|y\|^2\right)$ for all $x\in\R^n$.
\end{itemize}}
\end{remark}

\begin{remark}{\rm
The set $Y$ must be nonempty; otherwise, for every $x \in \R^n$ we have  by the convention $\sup\emptyset=-\infty$ that
	$$
		\sup\limits_{y\in Y} \left(J(x+y)-\dfrac{L}{2}\|y\|^2\right)+ \dfrac{L}{2}\inf\limits_{y\in Y} \|y\|^2\leq \sup\limits_{y\in Y} \left(J(x+y)-L\|y\|^2\right)=-\infty.
	$$
This implies that 
	$$	\inf_{x\in \R^n}\sup\limits_{y\in Y} \left(J(x+y)-\dfrac{L}{2}\|y\|^2\right)+ \dfrac{L}{2}\inf\limits_{y\in Y} \|y\|^2 \leq \inf_{x\in \R^n} J(x),$$
	even when the right-hand side  of the last inequality equals $-\infty$. As a result, condition~(c) cannot be fulfilled.}
\end{remark}

 Based on the above remarks, we can study Problem~1 under the assumption saying that $Y \subset \R^n$ is a nonempty closed convex set.

\medskip
The next minimax theorem for possibly nonconvex functions resolves Problem~1 in the negative. Namely, from the theorem it follows that one cannot find a $C^1$ function $J$ and a convex set $Y$ satisfying all the required conditions. 

\begin{theorem}\label{thm1}
Let $Y \subset \R^n$, with $n\geq 1$, be a nonempty closed convex set and let $J:\R^n \to \R$ be a Fr\'echet differentiable function satisfying the following conditions:
\begin{itemize}
	\item [\rm (a)] $\lim\limits_{\|x\|\to+\infty} J(x) =+\infty$;
	\item [\rm (b)] $\nabla J$ is Lipschitzian, with Lipschitz constant $L$.
\end{itemize}
Then, one has
\begin{equation}\label{eq1}
\inf_{x\in \R^n} J(x) -\dfrac{L}{2}\inf_{y\in Y}\|y\|^2= \inf_{x\in \R^n}\, \sup_{y\in Y} \left(J(x+y)-\dfrac{L}{2}\|y\|^2\right)
\end{equation} or, equivalently,
\begin{equation}\label{eq1-a}
\disp\sup_{y\in Y}\,\disp\inf_{x\in\mathbb R^n}{\mathcal{L}}(x,y)=\disp\inf_{x\in \mathbb R^n}\,\disp\sup_{y\in Y}{\mathcal{L}}(x,y)
\end{equation} 
with
\begin{equation}\label{L(x,y)}{\mathcal{L}}(x,y):= J(x+y)-\dfrac{L}{2}\|y\|^2\quad{\rm for}\ x,y\in\R^n.
\end{equation} 
Furthermore, the function ${\mathcal{L}}(x,y)$ has a saddle point on $\R^n\times Y$, that is, there exists $(\bar{x},\bar{y})\in \R^n\times Y$ satisfying 
\begin{equation}\label{eq2}
{\mathcal{L}}(\bar{x},y) \leq {\mathcal{L}}(\bar{x},\bar{y})\leq {\mathcal{L}}(x,\bar{y})
\end{equation}
for all $(x,y)\in \R^n\times Y$.
\end{theorem}
\noindent{\it Proof.} For the function ${\mathcal{L}}(x,y)$ defined by~\eqref{L(x,y)}, we  have 
\begin{equation*}\label{trans1}\begin{array}{rl}		
\disp\sup_{y\in Y}\, \inf_{x\in \R^n} {\mathcal{L}}(x,y) &=\disp\sup_{y\in Y}\, \disp\inf_{x\in \R^n} \left(J(x+y)-\dfrac{L}{2}\|y\|^2\right)\\ &=  \disp\sup_{y\in Y} \left(  \inf_{z\in \R^n} J(z)-\dfrac{L}{2}\|y\|^2\right) \\
&=  \disp\inf_{z\in \R^n} J(z)+\disp\sup_{y\in Y} \left(-\dfrac{L}{2}\|y\|^2\right) \\&= \disp\inf_{x\in \R^n} J(x) -\dfrac{L}{2}\inf_{y\in Y}\|y\|^2.
\end{array}
\end{equation*}
This proves the equivalence between~\eqref{eq1} and~\eqref{eq1-a}. In addition, since $$\sup\limits_{y\in Y}\, \inf\limits_{x\in \R^n} {\mathcal{L}}(x,y) \leq \inf\limits_{x\in \R^n}\, \sup\limits_{y\in Y} {\mathcal{L}}(x,y)$$ (see the remarks given in Section~\ref{Sect-1}), we have
\begin{equation}\label{eq3}
\inf_{x\in \R^n} J(x) -\dfrac{L}{2}\inf_{y\in Y}\|y\|^2 \leq \inf_{x\in \R^n}\, \sup_{y\in Y} {\mathcal{L}}(x,y) = \inf_{x\in \R^n}\, \sup_{y\in Y} \left(J(x+y)-\dfrac{L}{2}\|y\|^2\right).
\end{equation}
As the function $J:\R^n \to \R$ is continuous and coercive by the condition~(a), by applying the Weierstrass Theorem (see, e.g.,~\cite[Theorem~27.4]{Munkres_2000} or~\cite[Theorem~1, p.~15]{Polyak_1987}) to the restriction of $J$ on a suitable sublevel set, we can find a point $x_*\in\R^n$ such that 
$$J(x_*)= \min_{x\in \R^n} J(x).$$
Thus, by Fermat's rule (see~\cite[Theorem~1, p.~11]{Polyak_1987} and note that the result and its proof is also valid for a normed space setting) we obtain
\begin{equation}\label{eq4}
	\nabla J(x_*)=0.
\end{equation}
Let $\bar{y}$ be the metric projection of $0$ onto the nonempty closed convex set $Y$, i.e., $\bar{y}\in Y$ and $\|\bar{y}\|=\min\limits_{y\in Y}\|y\|.$ Then, by the projection theorem~\cite[Theorem~3.14]{BC_2011} we have
$$\langle 0 - \bar{y}, y-\bar{y}\rangle \leq 0\quad\mbox{\rm for all}\ y\in Y,$$
or, equivalently,
\begin{equation}\label{eq5}
\langle \bar{y},  y-\bar{y}\rangle \geq 0\quad\mbox{\rm for all}\ y\in Y.
\end{equation}
For $\bar{x}:= x_*-\bar{y}$, from~\eqref{L(x,y)} and~\eqref{eq4} we can deduce that
$$
\nabla_y (-{\mathcal{L}})(\bar{x},\bar{y})=L \bar{y} - \nabla J(\bar{x}+\bar{y})= L \bar{y}  - \nabla J(x_*)= L \bar{y},
$$
where $\nabla_y (-{\mathcal{L}})$ denotes the gradient of $(-{\mathcal{L}})$ with respect to the second variable. Combining this with~\eqref{eq5} yields
$$	\langle \nabla_y (-{\mathcal{L}})(\bar{x},\bar{y}), y-\bar{y}\rangle = \langle L \bar{y}, y-\bar{y}\rangle \geq 0,\quad\mbox{\rm for all}\, y\in Y.$$
Thus, one gets
\begin{equation}\label{optimalitycon}
0\in \nabla_y (-{\mathcal{L}})(\bar{x},\bar{y}) + N_Y(\bar{y}),
\end{equation}
where $N_Y(\bar{y}):=\{y^*\in\mathbb R^n\mid \langle y^*,y-\bar y\rangle\leq 0\}$ is the normal cone to $Y$ at $\bar{y}$ in the sense of convex analysis. 

\smallskip
\noindent{\sc Claim 1.} \textit{For every $x\in \R^n$, the function $(-{\mathcal{L}})(x,\cdot):\R^n \to \R$ is convex.}

\smallskip
Indeed, let $x\in \R^n$ be given arbitrarily. By~\eqref{L(x,y)} and the condition~(b), for any $y, y'\in \R^n$ one  has
\begin{equation*}
	\begin{array}{ll}
		&\langle \nabla_y (-{\mathcal{L}})(x,y) - \nabla_y (-{\mathcal{L}})(x,y'), y-y'\rangle \\&= \langle (L y - \nabla J(x+y)) - (L y' - \nabla J(x+y')), y-y'\rangle\\ 
		& = L\|y-y'\|^2 - \langle \nabla J(x+y)-\nabla J(x+y'),y-y'\rangle\\
		&\geq L\|y-y'\|^2  - \|\nabla J(x+y)-\nabla J(x+y')\|\|y-y'\|\\
		& \geq L\|y-y'\|^2  -  L\|(x+y)-(x-y')\| \|y-y'\|\\
		& = 0.
	\end{array}
\end{equation*}
This implies that the gradient mapping $\nabla_y (-{\mathcal{L}})(x,\cdot):\R^n \to \R^n$ is monotone. Hence, by~\cite[Proposition~17.10]{BC_2011} we conclude that the function $y\mapsto -{\mathcal{L}}(x,y)$ is convex.

\smallskip
Now, for every $x\in \R^n$, consider the convex optimization problem
\begin{equation}\label{px}
	\min\,\{-\mathcal{L}(x,y)\mid y\in Y\}. \tag{$P_x$}
\end{equation}
Since $\bar{y}$ satisfies the inclusion~\eqref{optimalitycon},  by the basic necessary and sufficient optimality conditions for optimal
solutions to a minimization problem with a convex objective function and a convex constraint (see, e.g.,~\cite[Theorem~7.15]{MN_2023}) we can assert that $\bar{y}$ is a solution of~\eqref{px} with $x=\bar{x}$. Thus,  $$-\mathcal{L}(\bar{x},\bar{y})= \min\limits_{y\in Y}(-\mathcal{L}(\bar{x},y)).$$ Consequently,
\begin{equation}\label{eq6}
\begin{array}{rcl}
\mathcal{L}(\bar{x},\bar{y})= \max\limits_{y\in Y}\,\mathcal{L}(\bar{x},y) &=& \max\limits_{y\in Y}\left(J(\bar{x}+y)-\dfrac{L}{2}\|y\|^2\right)\\& \geq & \inf\limits_{x\in \R^n}\, \sup\limits_{y\in Y} \left(J(x+y)-\dfrac{L}{2}\|y\|^2\right).
\end{array}
\end{equation}
Moreover, by the constructions of $x_*$ and $\bar x$ one has
$$\mathcal{L}(\bar{x},\bar{y})= J(\bar{x}+\bar{y})-\dfrac{L}{2}\|\bar{y}\|^2 = J(x_*)- \dfrac{L}{2}\|\bar{y}\|^2 = \inf_{x\in \R^n} J(x) -\dfrac{L}{2}\inf_{y\in Y}\|y\|^2.$$
Therefore, from~\eqref{eq6} it follows that
$$\inf_{x\in \R^n} J(x) -\dfrac{L}{2}\inf_{y\in Y}\|y\|^2\geq \inf_{x\in \R^n}\, \sup_{y\in Y} \left(J(x+y)-\dfrac{L}{2}\|y\|^2\right).$$
Together with~\eqref{eq3}, this gives the equality~\eqref{eq1}, which is equivalent to~\eqref{eq1-a}.

To verify the last assertion of the theorem, we observe for all $x \in \R^n$ that
\begin{align*}
{\mathcal{L}}(\bar{x},\bar{y}) = J(\bar{x}+\bar{y})-\dfrac{L}{2}\|\bar{y}\|^2& = J(x_*)-\dfrac{L}{2}\|\bar{y}\|^2 \\&\leq J(x+\bar{y})-\dfrac{L}{2}\|\bar{y}\|^2\\&={\mathcal{L}}(x,\bar{y}).
\end{align*}
Combining this with the fact that $\mathcal{L}(\bar{x},\bar{y})= \max\limits_{y\in Y}\,\mathcal{L}(\bar{x},y)$ establishes the inequalities in~\eqref{eq2} for all $(x,y)\in \R^n\times Y$. Thus, $(\bar{x},\bar{y})$ is a saddle point of the function ${\mathcal{L}}(x,y)$ on $\R^n\times Y$.
$\hfill\Box$

\medskip
Since all the arguments of the above proof remain valid if instead of $\mathbb{R}^n$ we consider a real Hilbert space $\cal{H}$ and replace the coercivity condition~(a) by the requirement saying that the optimization problem $\min\{J(x)\mid x\in \cal{H}\}$ has a solution $x_*$, the next minimax theorem for possibly nonconvex functions in a strong form is valid. 

\begin{theorem}\label{minimax_thm2} Let $\mathcal{H}$ be Hilbert space and $Y\subset \mathcal{H}$ be a nonempty closed convex set. Suppose that $J:\cal{H} \to \R$ is such a Fr\'echet differentiable function that $\nabla J$ is Lipschitzian, with Lipschitz constant $L$, and there is a point $x_*\in \mathcal{H}$ with $J(x_*)\leq J(x)$ for all $x\in \cal{H}$. Then, one has 
	\begin{equation*}
		\inf_{x\in \cal{H}} J(x) -\dfrac{L}{2}\inf_{y\in Y}\|y\|^2= \inf_{x\in \cal{H}}\, \sup_{y\in Y} \left(J(x+y)-\dfrac{L}{2}\|y\|^2\right),
	\end{equation*} which is equivalent to 
	\begin{equation*}\label{eq1b}
		\disp\sup_{y\in Y}\,\disp\inf_{x\in\cal{H}}{\mathcal{L}}(x,y)=\disp\inf_{x\in \cal{H}}\,\disp\sup_{y\in Y}{\mathcal{L}}(x,y)
	\end{equation*} with 
	$${\mathcal{L}}(x,y):= J(x+y)-\dfrac{L}{2}\|y\|^2\quad\mbox{\rm for}\  x,y\in\cal{H}.$$
	Moreover, the function ${\mathcal{L}}(x,y)$ possesses a saddle point, that is, there exists $(\bar{x},\bar{y})\in {\cal{H}}\times Y$ such that 
	\begin{equation*}
		{\mathcal{L}}(\bar{x},y) \leq {\mathcal{L}}(\bar{x},\bar{y})\leq {\mathcal{L}}(x,\bar{y})
	\end{equation*}
	for all $(x,y)\in {\cal{H}}\times Y$.
\end{theorem}

By introducing a perturbation function\footnote{This idea was proposed to the authors by Professor Biagio Ricceri during a seminar discussion in October 2024.} (see the function $\gamma$ below), we can have the following extension for Theorem~\ref{thm1}.

\begin{theorem}\label{minimax_thm3}
	Let $Y \subset \R^n$, with $n\geq 1$, be a nonempty closed convex set, $\gamma:\R^n\to\R$ a concave continuous function, and $J:\R^n \to \R$ a Fr\'echet differentiable function satisfying the following conditions:
	\begin{itemize}
		\item [\rm (a)] $\lim\limits_{\|x\|\to+\infty} J(x) =+\infty$;
		\item [\rm (b)] $\nabla J$ is Lipschitzian, with Lipschitz constant $L$.
	\end{itemize}
	Then, one has the  equality 
	\begin{equation}\label{eq7}
		\inf_{x\in \R^n} J(x) +\sup_{y\in Y}\left(-\dfrac{L}{2}\|y\|^2+\gamma(y)\right)= \inf_{x\in \R^n}\, \sup_{y\in Y} \left(J(x+y)-\dfrac{L}{2}\|y\|^2+\gamma(y)\right),
	\end{equation}
   which is equivalent to 
	\begin{equation}\label{eq1d}
		\disp\sup_{y\in Y}\,\disp\inf_{x\in\mathbb R^n}f(x,y)=\disp\inf_{x\in \mathbb R^n}\,\disp\sup_{y\in Y}f(x,y)
	\end{equation} 
	with
	\begin{equation}\label{f} f(x,y):= J(x+y)-\dfrac{L}{2}\|y\|^2+\gamma(y)\quad{\rm for}\ x,y\in\R^n.\end{equation} 
	In addition, $f(x,y)$ admits a saddle point, that is, there exists a point $(\bar{x},\bar{y})\in \R^n\times Y$ such that 
	\begin{equation}\label{eq8}
		f(\bar{x},y) \leq f(\bar{x},\bar{y})\leq f(x,\bar{y})
	\end{equation}
	for all $(x,y)\in \R^n\times Y$.
\end{theorem}

\noindent{\it Proof.} For the function $f(x,y)$ defined by~\eqref{f}, note that
\begin{align*}
	\sup_{y\in Y}\, \inf_{x\in \R^n} f(x,y) &=\sup_{y\in Y}\, \inf_{x\in \R^n} \left(J(x+y)-\dfrac{L}{2}\|y\|^2+\gamma(y)\right)\\&=  \sup_{y\in Y} \left(\inf_{z\in \R^n} J(z) -\dfrac{L}{2}\|y\|^2+\gamma(y)\right) \\&= \inf_{x\in \R^n} J(x) +\sup_{y\in Y}\left(-\dfrac{L}{2}\|y\|^2+\gamma(y)\right).
\end{align*} So,~\eqref{eq7} is equivalent to~\eqref{eq1d}. Besides, since $\sup\limits_{y\in Y}\, \inf\limits_{x\in \R^n} f(x,y) \leq \inf\limits_{x\in \R^n}\, \sup\limits_{y\in Y} f(x,y)$ (see the remarks given in Section~\ref{Sect-1}), we obtain
\begin{equation}\label{eq9}
	\begin{array}{rcl}
		\inf\limits_{x\in \R^n} J(x) +\sup\limits_{y\in Y}\left(-\dfrac{L}{2}\|y\|^2+ \gamma(y)\right)& \leq& \inf\limits_{x\in \R^n}\, \sup\limits_{y\in Y} f(x,y)\\& =&\inf\limits_{x\in \R^n}\, \sup\limits_{y\in Y} \left(J(x+y)-\dfrac{L}{2}\|y\|^2+\gamma(y)\right).
	\end{array}
\end{equation}
Thanks to the continuity of $J:\R^n \to \R$ and the coercivity condition~(a), arguing as in the proof of Theorem~\ref{thm1}, we can find a point $x_*\in\R^n$  such that 
$$J(x_*)= \min_{x\in \R^n} J(x).$$
Then, applying Fermat's rule (see~\cite[Theorem~1, p.~11]{Polyak_1987}) we get
\begin{equation}\label{eq10}
	\nabla J(x_*)=0.
\end{equation}
Consider the following optimization problem
\begin{equation}\label{op_y}
	\min \;\left\{\dfrac{L}{2}\|y\|^2-\gamma(y) \mid y\in Y\right\}.
\end{equation}
Since $\gamma:\R^n\to \R$ is concave, $-\gamma(\cdot)$ is convex, and hence by~\cite[Proposition 6.2]{BC_2011} the function $y\mapsto \dfrac{L}{2}\|y\|^2-\gamma(y)$ is strongly convex. Thus, by~\cite[Lemma~6]{Nesterov_2009} one can find a unique solution $\bar{y}$ to the optimization problem~\eqref{op_y}. Therefore, it follows from the Fermat stationary rule for convex extended-real-valued functions (see, e.g.,~\cite[Proposition~3.29]{MN_2022}) that
$$0\in \partial \left(\dfrac{L}{2}\|\cdot\|^2-\gamma(\cdot) +\delta(\cdot;Y)\right)(\bar{y}),$$
where $\delta(y;Y)=0$ for $y\in Y$ and $\delta(y;Y)=+\infty$ for $y\notin Y$ is the \textit{indicator function} of $Y$, and $\partial \left(\dfrac{L}{2}\|\cdot\|^2 -\gamma(\cdot)+\delta(\cdot;Y)\right)(\bar{y})$ is the subdifferential of the convex function $$y\mapsto \dfrac{L}{2}\|y\|^2-\gamma(y)+\delta(y;Y)\quad\; (y\in\mathbb R^n)$$ at the point $\bar{y}$. Then, noting that both functions $-\gamma(\cdot)$ and $\dfrac{L}{2}\|\cdot\|^2$ are continuous on $\R^n$, by the Moreau-Rockafellar theorem (see, e.g.,~\cite[Theorem~3.48]{MN_2022} or~\cite[Theorem~3.18]{MN_2023}) we have 
\begin{equation}\label{eq11}
	0\in L \bar{y}+ \partial (-\gamma)(\bar{y})+N_Y(\bar{y}).
\end{equation}
Since  $- J(x, \cdot)+\dfrac{L}{2}\|\cdot\|^2 $ is convex  for all $x\in \R^n$ (see Claim~1 in the proof of Theorem~\ref{thm1}), the function $- J(x, \cdot)+\dfrac{L}{2}\|\cdot\|^2 -\gamma(\cdot)$ is also convex for all $x\in\R^n$. That is, for every $x\in\mathbb R^n$, the function $-f(x,\cdot)$ is convex on $\R^n$. Put $\bar{x}= x_*-\bar{y}$. By the Moreau-Rockafellar theorem (see, e.g.,~\cite[Theorem~3.48]{MN_2022} or~\cite[Theorem~3.18]{MN_2023}),~\eqref{f}, and~\eqref{eq10}, we have
$$
\begin{array}{ll}
\partial_y (-f)(\bar{x},\bar{y})=- \nabla J(\bar{x}+\bar{y})  +L \bar{y}+ \partial (-\gamma)(\bar{y})&=- \nabla J(x_*)  + L \bar{y}  +  \partial (-\gamma)(\bar{y})\\&= L \bar{y}+\partial (-\gamma)(\bar{y}),
\end{array}
$$
where $\partial_y (-f)(\bar{x},\bar{y})$ denotes the subdifferential of $(-f)(\bar x,.)$ at $\bar{y}$. Combining this with~\eqref{eq11} yields
\begin{equation}\label{op_con}
	0\in \partial_y (-f)(\bar{x},\bar{y})+ N_Y(\bar{y}).
\end{equation}

\smallskip
Now, for every $x\in \R^n$, consider the convex optimization problem
\begin{equation}\label{op_prob}
	\min\,\{-f(x,y)\mid y\in Y\}. \tag{$P^\prime_x$}
\end{equation}
 Since $\bar{y}$ satisfies~\eqref{op_con}, it follows from the basic necessary and sufficient optimality conditions for optimal
 solutions to a minimization problem with a convex objective function and a convex constraint (see, e.g.,~\cite[Theorem~7.15]{MN_2023}) that  $\bar{y}$ is a solution of~\eqref{op_prob} with $x=\bar{x}$. Thus,  $$-f(\bar{x},\bar{y})= \min\limits_{y\in Y}(-f(\bar{x},y)).$$ Therefore,
\begin{equation}\label{eq12}
	\begin{array}{rcl}
		f(\bar{x},\bar{y})= \max\limits_{y\in Y}\,f(\bar{x},y) &=& \max\limits_{y\in Y}\left(J(\bar{x}+y)-\dfrac{L}{2}\|y\|^2+\gamma(y)\right)\\& \geq & \inf\limits_{x\in \R^n}\left(\sup\limits_{y\in Y} \left(J(x+y)-\dfrac{L}{2}\|y\|^2+\gamma(y)\right)\right) .
	\end{array}
\end{equation}
In addition, as $\bar{y}$ is the solution of~\eqref{op_y}, we have
\begin{equation*}
\begin{array}{rcl}
f(\bar{x},\bar{y})= J(\bar{x}+\bar{y})-\dfrac{L}{2}\|\bar{y}\|^2+\gamma(\bar{y}) &=& J(x_*)- \dfrac{L}{2}\|\bar{y}\|^2 +\gamma(\bar{y})\\&= &\inf\limits_{x\in \R^n} J(x) +\sup\limits_{y\in Y} \left(-\dfrac{L}{2}\|y\|^2+\gamma(y)\right).
\end{array}
\end{equation*}
By~\eqref{eq12}, it holds that
$$\inf\limits_{x\in \R^n} J(x) +\sup\limits_{y\in Y} \left(-\dfrac{L}{2}\|y\|^2+\gamma(y)\right)\geq \inf\limits_{x\in \R^n}\, \sup\limits_{y\in Y} \left(J(x+y)-\dfrac{L}{2}\|y\|^2+\gamma(y)\right).$$
This in combination with~\eqref{eq9} gives~\eqref{eq7}.

It remains to verify the inequalities in~\eqref{eq8} for every $(x,y)\in \R^n\times Y$. Note that for every $x \in \R^n$ we have
\begin{align*}
	f(\bar{x},\bar{y}) = J(\bar{x}+\bar{y})-\dfrac{L}{2}\|\bar{y}\|^2+\gamma(\bar{y})& = J(x_*)-\dfrac{L}{2}\|\bar{y}\|^2 +\gamma(\bar{y})\\&\leq J(x+\bar{y})-\dfrac{L}{2}\|\bar{y}\|^2+\gamma(\bar{y})\\&=f(x,\bar{y}).
\end{align*}
Together with the fact that $f(\bar{x},\bar{y})= \max\limits_{y\in Y}\, f(\bar{x},y)$, this implies~\eqref{eq8}.
$\hfill\Box$

\section{An Illustrative Example}\label{Sect-3}

As an illustration for the three minimax theorems for possibly nonconvex functions, which were obtained in Section~\ref{Sect-2}, we now construct a concrete example. Since Theorem~\ref{thm1} is a special case of Theorems~\ref{minimax_thm2} and~\ref{minimax_thm3}, it suffices to refer to Theorem~\ref{thm1}.

\medskip
In the following example, $J$ is an one-variable $C^1$ function, which is neither convex nor concave. Since some sublevel sets of the function is nonconvex, it is not quasiconvex. Moreover, as two stationary points of the function are not global minimizers, it is not a pseudoconvex function.  

\begin{example}\label{example1}
In the setting of Theorem~\ref{thm1}, let $n=1$, $Y=\R$, and the function $J:\R\to\R$ be defined by
\begin{equation}\label{example1_function}
J(x)=
\begin{cases}
(x+\pi)^2 - (x+\pi)\qquad &\text{if}\;\; x\leq -\pi\\
\sin x &\text{if}\;\; -\pi< x<\pi\\
(x-\pi)^2-(x-\pi)&\text{if}\;\; x\geq \pi
\end{cases}
\end{equation}
\begin{figure}[h]
	\centering
	\includegraphics[width=0.8\textwidth]{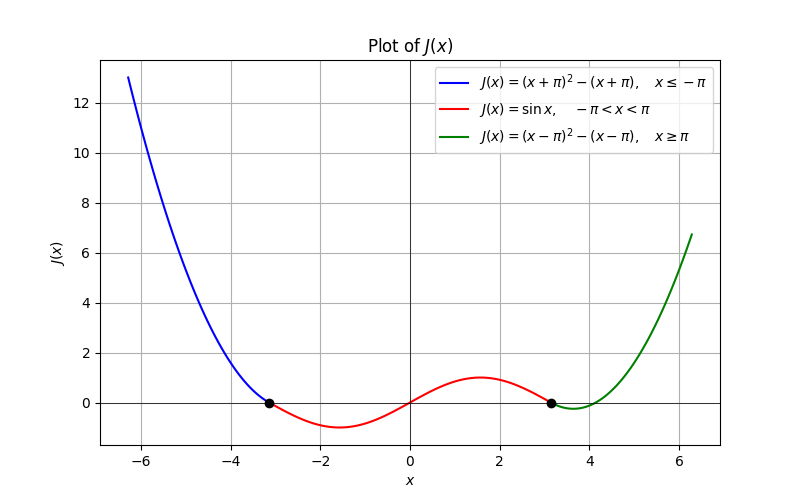}
	\caption{Graph of $J$ defined in Example~\ref{example1}}
	\label{fig:example1}
\end{figure}
(see Fig.~1). Then, $J$ is a nonconvex function and we have
\begin{equation}\label{example1_grad}
	\nabla J(x)=
	\begin{cases}
		2(x+\pi)-1\qquad &\text{if}\;\; x\leq -\pi\\
		\cos x &\text{if}\;\;  -\pi< x< \pi\\
		2(x-\pi)-1&\text{if}\;\; x\geq \pi.
	\end{cases}
\end{equation}
From~\eqref{example1_grad} it follows that $J$ is a $C^1$ function. By~\eqref{example1_function}, it is obvious that
$\lim\limits_{|x|\to+\infty} J(x) =+\infty;$ hence $J$ is a coercive function.

\smallskip 
\noindent {\sc Claim 1.} \textit{$\nabla J$ is Lipschitz continuous, with Lipschitz constant $L=2$.} 

\smallskip
Indeed, given two district points $x,y\in \R$, we can assume without any loss of generality that $x<y$. When $x$ and $y$ belong to the same interval, either $(-\infty,-\pi]$, or $(-\pi,\pi)$, or $[\pi,+\infty)$, by~\eqref{example1_grad} we can infer that
\begin{equation}\label{ex1_case0}
	|\nabla J(x)-\nabla J(y)|  \leq 2(y-x).
\end{equation} Otherwise, the intersection $\Omega:=(x,y)\cap \{-\pi, \pi\}$ must be nonempty. In the case $\Omega=\{-\pi\}$, by~\eqref{example1_grad} we have 
\begin{equation*}\begin{aligned}
	|\nabla J(x)-\nabla J(y)|&=|\big(\nabla J(x)-\nabla J(-\pi)\big)+\big(\nabla J(-\pi)-\nabla J(y)\big)|\\
	&\leq |\nabla J(x)-\nabla J(-\pi)|+|\nabla J(-\pi)-\nabla J(y)|\\
	& \leq 2\big((-\pi)-x\big)+2\big(y-(-\pi)\big)\\
	&= 2(y-x).
\end{aligned}
\end{equation*} If $\Omega=\{-\pi\}$, or $\Omega=\{-\pi, \pi\}$, then the inequality~\eqref{ex1_case0} can be proved analogously. So, $\nabla J$ is Lipschitz continuous, with Lipschitz constant $L=2$.

\smallskip
We have shown that the nonconvex function $J$ satisfies all the assumptions stated in Theorem~\ref{thm1}.

\smallskip
The function ${\mathcal{L}}(x,y)$ in~\eqref{L(x,y)} now has the form
\begin{equation}\label{L_exam}{\mathcal{L}}(x,y)= J(x+y)-y^2\quad{\rm for}\ x,y\in\R,
\end{equation} because $n=1$ and $L=2$. By definition, $(\bar x,\bar y)\in \R\times\R$ is a saddle point of ${\mathcal{L}}(x,y)$ on $\R\times Y=\R\times\R$ if and only if
\begin{equation}\label{eq2x}
{\mathcal{L}}(\bar{x},\bar{y})\leq {\mathcal{L}}(x,\bar{y})\quad\; \forall x\in\R
\end{equation}
and
\begin{equation}\label{eq2y}
	{\mathcal{L}}(\bar{x},y) \leq {\mathcal{L}}(\bar{x},\bar{y})\quad\; \forall y\in\R.
\end{equation} Let $\bar z:=\bar x+\bar y$. 

\smallskip
Using~\eqref{L_exam}, we see that the condition~\eqref{eq2x} is satisfied whenever~$\bar x$ is the global minimizer of the function $x\mapsto J(x+\bar y)$ on $\R$. Clearly, the latter means that $\bar z$ is the global minimizer of the function $z\mapsto J(z)$ on $\R$.  So, by~\eqref{example1_function} (see also the graph of $J$ in Fig.~1) we can assert that~\eqref{eq2x} is fulfilled if and only if $\bar z=-\dfrac{\pi}{2}$.

\smallskip
Next, note that the condition~\eqref{eq2y} is satisfied whenever~$\bar y$ is the global minimizer of the function $(-{\mathcal{L}})(\bar x,\cdot):\R^n \to \R$. By Claim~1 in the proof of Theorem~\ref{thm1}, we know that $(-{\mathcal{L}})(\bar x,\cdot)$ is a convex function. Hence, by Fermat's rule, the fulfillment of~\eqref{eq2y} is equivalent to the condition $\nabla_y(-{\mathcal{L}})(\bar x,\bar y)=0$. In accordance with~\eqref{L_exam}, the latter can be rewritten as
\begin{equation}\label{stationary_y}
	-\nabla J(\bar x+\bar y)+2\bar y=0.
\end{equation} Since $\bar x+\bar y=\bar z=-\frac{\pi}{2}$, from~\eqref {example1_function} we can deduce that $\nabla J(\bar x+\bar y)=0$. Therefore,~\eqref{stationary_y} means that $\bar y=0$. Substituting this value $\bar y$ into the equality $\bar x+\bar y=-\frac{\pi}{2}$ yields $\bar x=-\frac{\pi}{2}$.

\smallskip
Summing up all the above, we conclude that the function ${\mathcal{L}}(x,y)$ in~\eqref{L_exam} has the unique saddle point $(\bar x,\bar y)=\left(-\frac{\pi}{2},0\right)$ on $\R\times Y=\R\times\R$, with the saddle value ${\mathcal{L}}(\bar x,\bar y)=-1$. As a consequence, the equality~\eqref{eq1-a} holds and one has \begin{equation*}\label{eq1a}
	\disp\sup_{y\in Y}\,\disp\inf_{x\in\mathbb R}{\mathcal{L}}(x,y)=\disp\inf_{x\in \mathbb R}\,\disp\sup_{y\in Y}{\mathcal{L}}(x,y)=-1
\end{equation*} 
(see the remarks given in Section~\ref{Sect-1}). So, just by using direct computation, we have verified all the assertions of Theorem~\ref{thm1} for the function $J$ given by~\eqref{example1_function} and the set $Y=\R$. 
\end{example}

\section{Concluding Remarks}\label{Sect-4}

Three minimax theorems for possibly nonconvex functions on $\mathbb R^n$ or on infinite-dimensional Hilbert spaces have been obtained in this paper. The first theorem provides a comprehensive solution to an interesting open problem related to continuously differentiable functions. 

\medskip
A carefully designed example illustrates the given minimax theorems.

\medskip
Extensions of Theorem~\ref{thm1}--\ref{minimax_thm3} to the case where $J$ need not be a $C^1$ function (say, $J$ is just a locally Lipschitz function) deserve further investigations.
\bigskip

\noindent{\bf Acknowledgements:} This research was supported by the project NCXS02.01/24-25 of the Vietnam Academy of Science and Technology. The authors would like to thank Professor Biagio Ricceri for proposing the open problem studied in the present paper and for warm hospitality at Catania University, Italy.

\end{document}